\documentclass[graybox]{svmult}

\usepackage{mathptmx}       
\usepackage{helvet}         
\usepackage{courier}        
\usepackage{type1cm}        
\usepackage{makeidx}         
\usepackage{graphicx}        
\usepackage{multicol}        
\usepackage[bottom]{footmisc}

\usepackage{amssymb,amsbsy,amsmath}
\usepackage{afterpage}

\makeindex             

\begin{document}

\title*{Substitution rules and topological properties of the Robinson tilings}
\author{Franz G\"ahler}
\institute{Franz G\"ahler \at Faculty of Mathematics, University of Bielefeld, 
33615 Bielefeld, Germany \\ \email{gaehler@math.uni-bielefeld.de}}
%
%
\maketitle

\abstract{A relatively simple substitution for the Robinson tilings is
presented, which requires only 56 tiles up to translation. In this
substitution, due to Joan M.\ Taylor, neighboring tiles are substituted
by partially overlapping patches of tiles. We show that this overlapping 
substitution gives rise to a normal primitive substitution as well,
implying that the Robinson tilings form a model set and thus have pure 
point diffraction. This substitution is used to compute the 
\v{C}ech cohomology of the hull of the Robinson tilings via the 
Anderson-Putnam method, and also the dynamical zeta function of the 
substitution action on the hull. The dynamical zeta function is then 
used to obtain a detailed description of the structure of the hull,
relating it to features of the cohomology groups.
}


\section{Introduction}
\label{sec:intro}

Robinson's aperiodic set of tiles \cite{Rob71} was the first reasonably 
small such set which could tile the plane only aperiodically. The 
local matching rules enforce a hierarchical structure into the tilings, 
which is used to prove that only aperiodic tilings are admitted. Despite
this hierarchical structure, for a long time it was not known whether
the Robinson tilings can be generated also by a substitution, which
would have enormous advantages for a more detailed study. Only very 
recently, a substitution for the Robinson tilings could be constructed
explicitly \cite{GJS12}, albeit a rather complicated one. The Robinson 
tilings therefore remain an interesting example, not only for historical 
reasons. In this paper, we present a much simpler substitution, derived 
from an overlapping substitution due to Joan M.~Taylor, which we then use 
to analyse the structure of the hull of the Robinson tilings in more 
detail, and relate it to some of the topological invariants of the hull.


\section{A simple substitution for the Robinson tilings}
\label{sec:subst}

\begin{figure}[t]
\includegraphics[height=1cm]{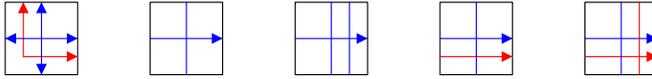}
\caption{The Robinson tiles. The one on the left is called a cross
and plays a special role. All tiles can also be rotated and reflected.}
\label{fig:robtiles}
\end{figure}

\begin{figure}[t]
\sidecaption
\includegraphics[width=6.5cm]{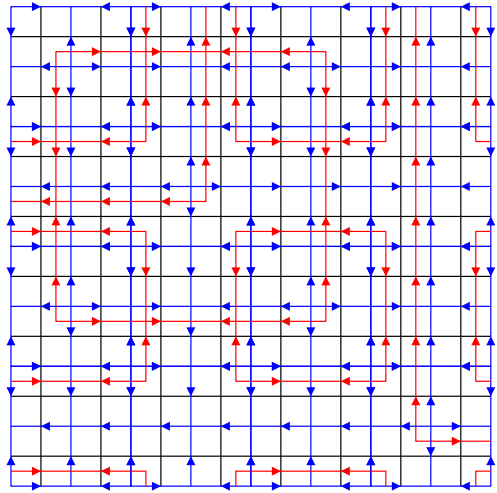}
\caption{A patch of a Robinson tiling. Note the red square frames with 
corners at cross tiles, which occur at all sizes of the form $2^n$, 
proving the aperiodicity of the Robinson tilings. The corners of the 
smallest square frames are at crosses forming the odd/odd sublattice 
of tiles.}
\label{fig:robtiling}
\end{figure}

Robinson tilings consist of the five square tiles shown in
Fig.~\ref{fig:robtiles}. As the tiles are allowed to be rotated and 
reflected, there are 28 tiles up to translation. In a legal Robinson 
tiling, the tiles must obey some local rules. Firstly, the decoration 
lines must continue across edges, with exactly one arrow head at each 
line join. Secondly, there must be a square sublattice of index 
4 whose tiles are all cross tiles. Apart from this lattice of cross tiles,
there may be other crosses as well. We assume in the following, that 
this sublattice of cross tiles is a the odd/odd position. All tilings 
satisfying the two rules (which are both local) are called Robinson 
tilings. In any Robinson tiling, the decoration lines form a hierarchy 
of square frames of all sizes $2^n$ (see Fig.~\ref{fig:robtiling}), 
which proves that Robinson tilings cannot be periodic.

The local rules given above admit also some tilings with defect
lines, which are not repetitive (for details, see \cite{Rob71,JM97}). 
As we are heading at primitive substitution rules, by which we can 
reach only repetitive tilings, we want to discard these defective 
tilings. We therefore confine ourselves to the minimal subspace of 
repetitive tilings which is closed and invariant under translations 
and substitutions. The tilings which we discard form a set of measure 
zero. In particular, their exclusion does not change any spectral 
properties. 

The hierarchy of square frames of all sizes (Fig.~\ref{fig:robtiling})
suggests a hierarchical structure in the tilings, and it would only
be natural if the Robinson tilings could be constructed also by a
substitution rule. The construction of such a substitution rule was
achieved only recently \cite{GJS12}, and the substitution proved to be 
rather complicated, with 208 tiles up to translation. The reason is,
that the self-similarity inherent in the Robinson tilings scales
around the tile centers, not the vertices. For the substitution,
one therefore had to disect and reassemble the original tiles to
new ones, having their vertices at the original tile centers, which
results in the rather large number of tiles.

\begin{figure}[P]
\centerline{
  \includegraphics[width=5.2cm]{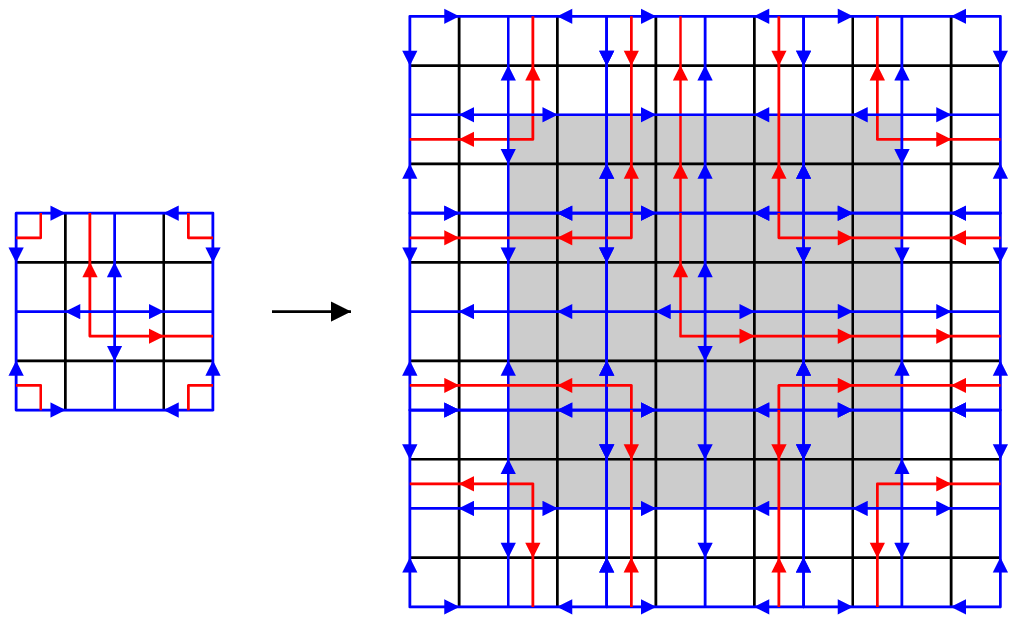}\hfill
  \includegraphics[width=5.2cm]{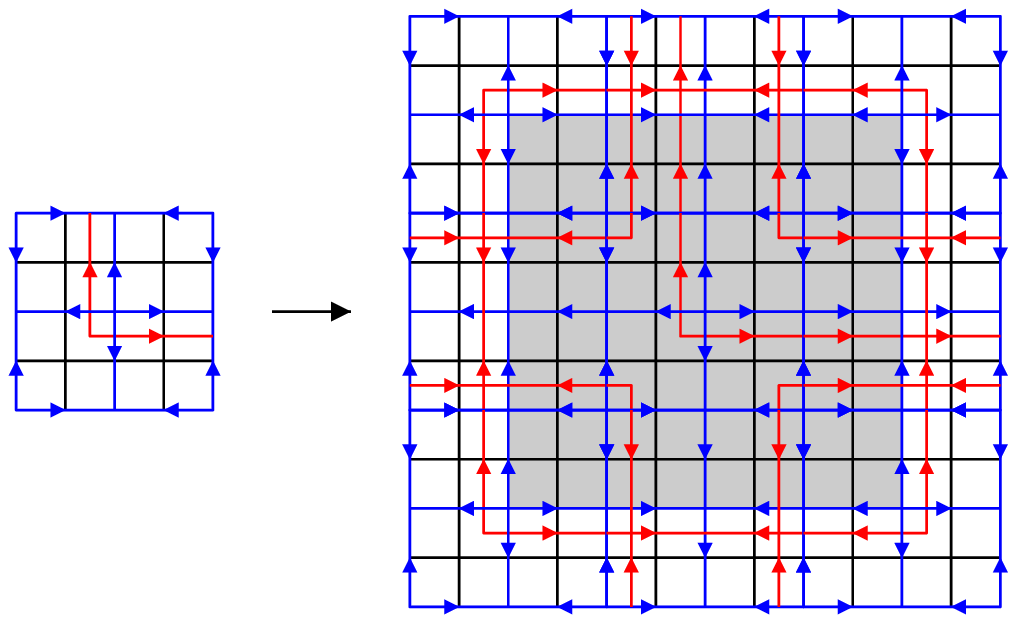}
}
\medskip

\centerline{
  \includegraphics[width=5.2cm]{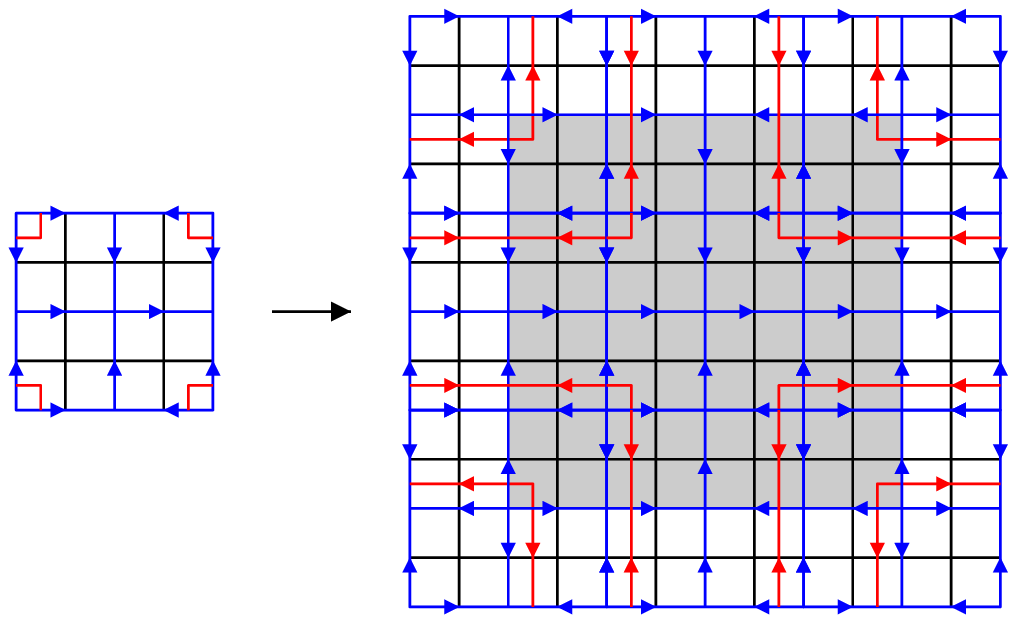}\hfill
  \includegraphics[width=5.2cm]{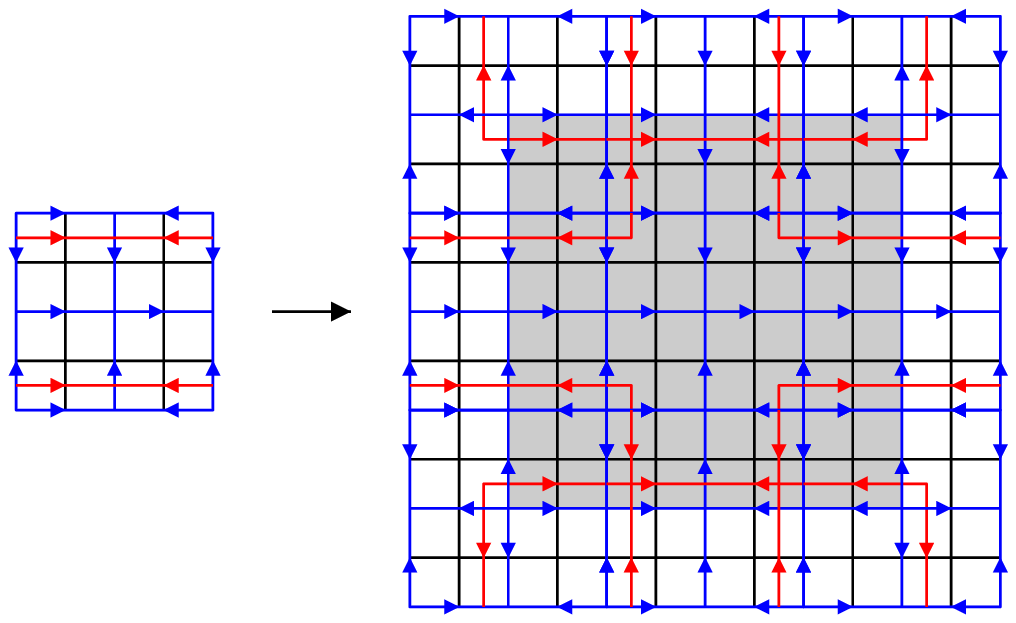}
}
\medskip

\centerline{
  \includegraphics[width=5.2cm]{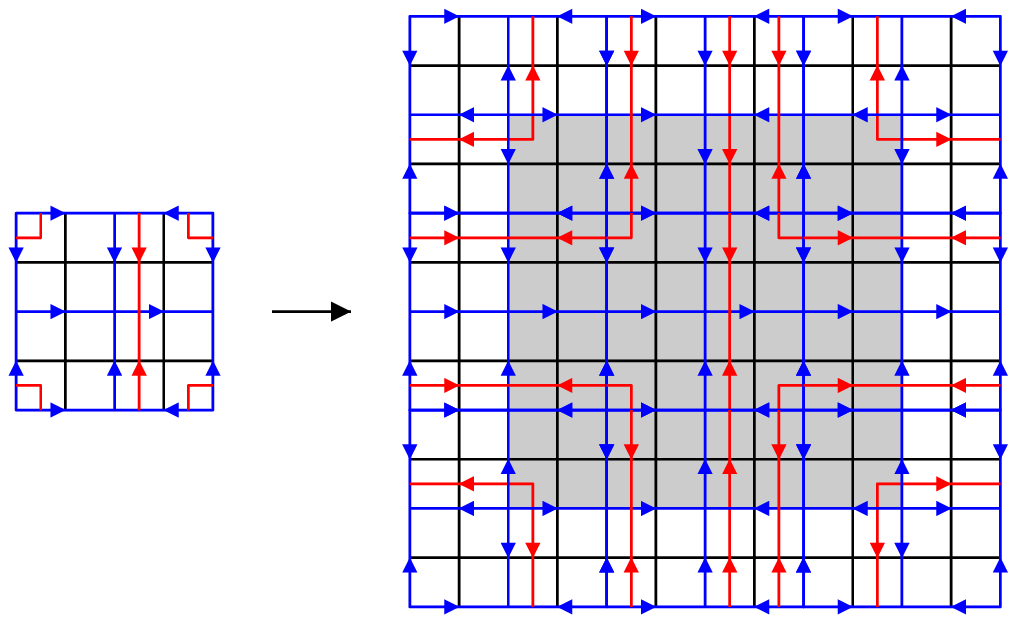}\hfill
  \includegraphics[width=5.2cm]{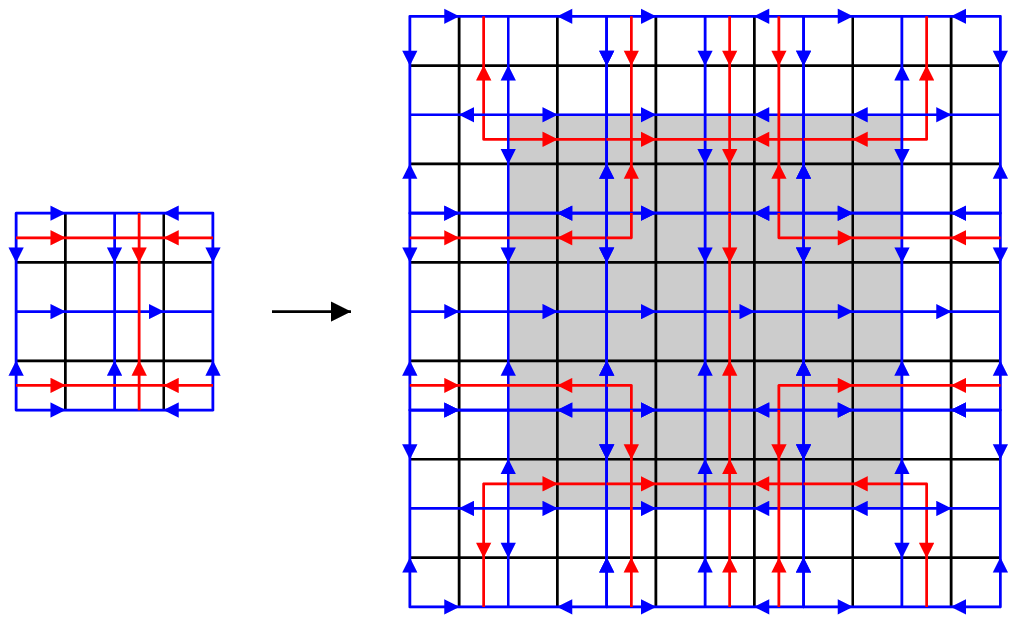}
}
\medskip

\centerline{
  \includegraphics[width=5.2cm]{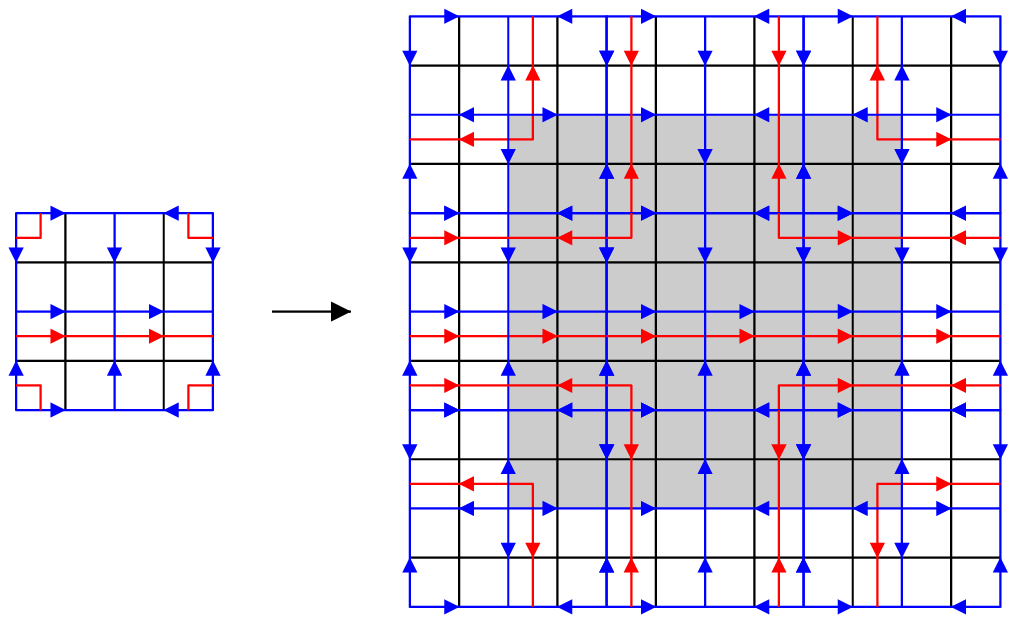}\hfill
  \includegraphics[width=5.2cm]{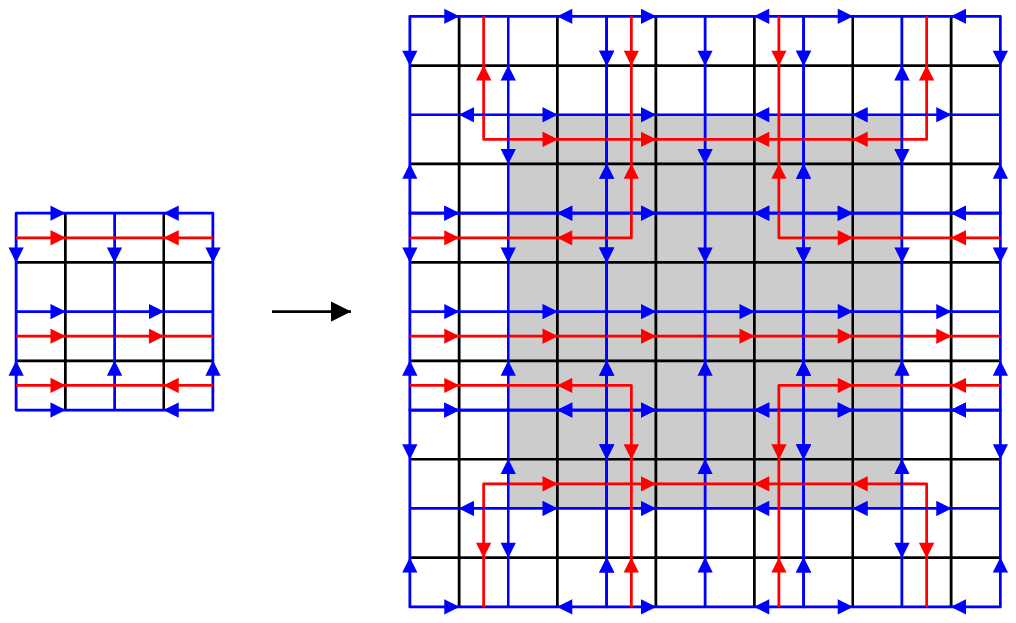}
}
\medskip

\centerline{
  \includegraphics[width=5.2cm]{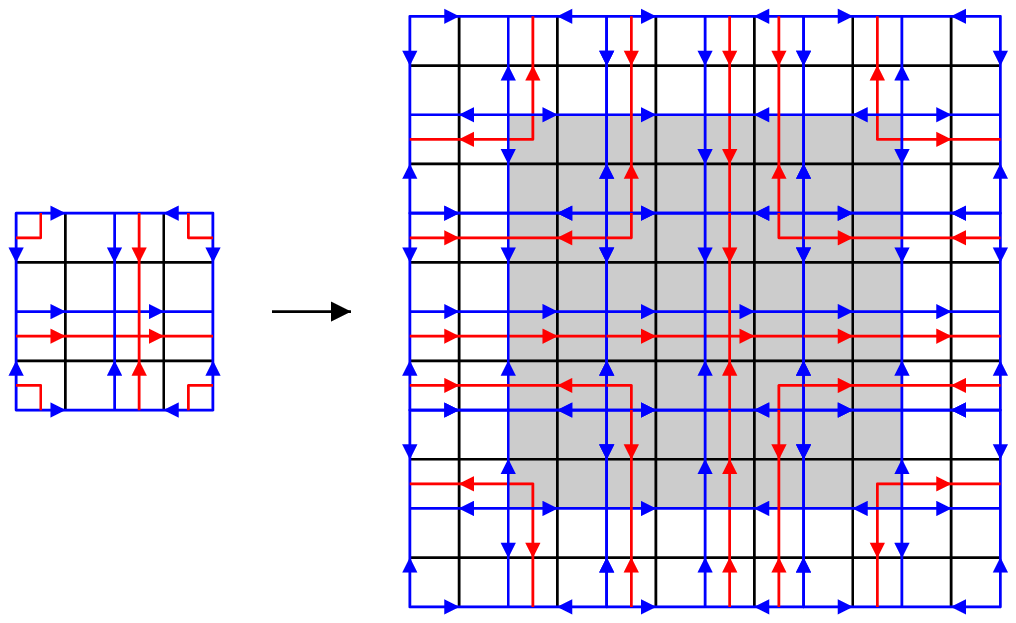}\hfill
  \includegraphics[width=5.2cm]{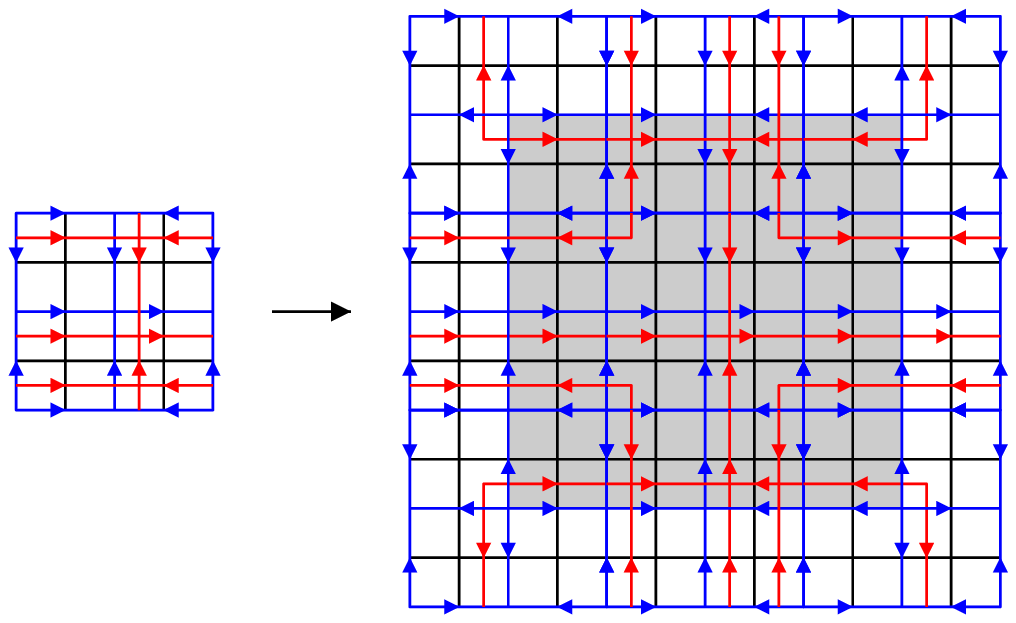}
}
\caption{Overlapping substitution for the Robinson tiling.  Each tile
  is replaced by a $3\times3$-patch of tiles.  Rotated/reflected
  tiles are substituted by the corresponding rotated/reflected patches.
  The inflated tiles cover only the area shaded in gray.  The 
  substitutions of neighboring tiles have thus an overlap, on which 
  they agree.}
\label{fig:robsubst}
\end{figure}

Here, we want to follow a different route, starting from a proposal
of Joan M.~Taylor (private communication). Recall that the self-similarity
scales about tile centers. The idea now is to replace a tile by a
$3\times3$-patch of tiles under the substitution. This patch is larger 
than the original tile inflated by a factor of 2, so that there are 
consistency conditions to be obeyed: the substitutions of neighboring 
tiles have an overlap, on which they must agree. A relatively simple 
solution is obtained if we pass to new tiles which are larger by a 
factor 2. These new tiles have their centers at the tiles at even/even 
positions (recall that the tiles at odd/odd positions are all crosses).
If we add to those even/even tiles a layer of thickness one half,
all the remaining tiles are consumed, and we end up with new square
tiles of edge length 2 at even/even positions. It turns out that the 
28 translation classes of tiles at even/even positions split up into 
two classes each, so that we now have 56 tile types up to translation. 
Moreover, these tiles admit a well-defined overlapping substitution, 
as shown in Fig.~\ref{fig:robsubst}.

\afterpage{\clearpage}

The overlapping substitution of Fig.~\ref{fig:robsubst} is considerably
simpler than the one found previously \cite{GJS12}. The set of translation
classes of tiles has been cut to a mere 56, from 208 previously.
For certain applications, however, such as the computation of the 
cohomology via the Anderson-Putnam method \cite{AP98}, an overlapping 
substitution is not suitable. To avoid this problem, we observe
that we can always pass to a normal substitution by replacing a tile
not by a full $3\times3$-patch, but by the $2\times2$-subpatch at the
upper right corner, say. Note that we always have to take the subpatch 
at the same corner, also for the rotated tiles, so that each tile is 
assigned to a unique supertile. As a result, this assignment breaks 
the rotation/reflection covariance of the substitution rules, but this 
is a small price to pay. 

Having derived our substitution from an overlapping substitution has yet 
another advantage. Since the $3\times3$-patches cover more than the 
inflated tiles, the overlapping substitution obviously forces the border 
\cite{Kel95}, a property which is inherited also by the normal substitution
derived from it. This allows to avoid the use of collared tiles \cite{AP98}
in the Anderson-Putnam method, which is a tremendous advantage, as it
also helps to keep the number of tile types small.


\section{The structure of the hull}
\label{sec:struct}

Due to the repetitivity, the translation group acts minimally on
the space of all repetitive Robinson tilings: every translation
orbit is dense. The tiling space is therefore the hull of any of
its member tilings. Having a substitution, the hull can now be 
constructed as an inverse limit space \cite{AP98}, and having a 
{\em simple} substitution which forces the border and requires 
only 56 tiles up to translation simplifies the task considerably. 
The mere fact of having a lattice substitution tiling has some 
immediate consequences. Since the crosses at odd/odd positions 
form a lattice-periodic subset of tiles (with period 4 in each
direction), results of Lee and Moody 
\cite{LM01} allow to conclude that the Robinson tilings form a 
model set and are thus pure-point diffractive. Since the defective 
Robinson tilings are a subset of measure zero, the pure-point 
diffractiveness extends even to all Robinson tilings.

As a limit-periodic model set, the space of Robinson tilings must
project 1-to-1 almost everywhere to an underlying 2d, 2-adic solenoid
$\mathbb{S}^2_2$ via the torus parametrisation \cite{BLM07}.
In the following, we will analyse the structure of the set where 
this projection fails to be 1-to-1, and try to connect it to the
\v{C}ech cohomology of the hull. The latter was obtained
in \cite{GJS12} via the Anderson-Putnam method \cite{AP98} as
\begin{equation}
H^2=\mathbb{Z}[{\scriptstyle\frac14}]
    \oplus\mathbb{Z}[{\scriptstyle\frac12}]^{10}
    \oplus\mathbb{Z}^8\oplus\mathbb{Z}_4,\quad 
     H^1=\mathbb{Z}[{\scriptstyle\frac12}]^2\oplus\mathbb{Z},
     \quad H^0=\mathbb{Z},
\end{equation}
which is confirmed using our new, simpler substitution. There
is a natural substitution action on the hull, whose Artin-Mazur
zeta function is defined as
\begin{equation}
  \zeta(z) = \exp\left(\sum_{m=1}^\infty \frac{a_m}{m} z^m\right) 
\end{equation}
where $a_m$ is the number of points in the hull that are invariant 
under an $m$-fold substitution. Note that if the hull consists of
two components for which the periodic points can be counted separately,
$a_m=a'_m+a''_m$, the corresponding partial zeta functions have to be 
multiplied: $\zeta(z) = \zeta'(z)\cdot\zeta''(z)$.

Anderson and Putnam have given a different way to compute the dynamical
zeta function, as a by-product of computing the \v{C}ech cohomology
\cite{AP98}. Recall that the hull is obtained as the inverse limit
of the substitution acting on an approximant cell complex. As a
consequence, the cohomology of the hull is the direct limit of
the substitution action on the cohomology of that cell complex.
Suppose $A^{(m)}$ is the matrix of the substitution action on
the $m$-th cohomology group (with rational coefficients) of the
hull of a substitution tiling. The dynamical zeta function is then
given by \cite{AP98}
\begin{equation}
 \zeta(z) = \frac{\prod_{k\ \text{odd}}\det(1-zA^{d-k})}
                   {\prod_{k\ \text{even}}\det(1-zA^{d-k})}
          = \frac{\prod_{k\ \text{odd}}\prod_i(1-z\lambda_i^{d-k})}
                   {\prod_{k\ \text{even}}\prod_i(1-z\lambda_i^{d-k})}
\end{equation}
where the latter equality holds if the matrices $A^{(m)}$ diagonalizable,
and the $\lambda_i^{(m)}$ are their eigenvalues. Note that Anderson
and Putnam have used the matrices of the substitution action on
the cochain groups of the approximant complex, rather than the cohomology,
but the additional terms in their formula cancel between numerator and 
denominator.

If we apply this to the Robinson tilings, and take into account
the eigenvalues of the substitution action on the cohomology,
we obtain for the zeta function
\begin{align} 
  \zeta(z) &= \frac{ (1-2z)^2 (1-z) }
      { (1-z) (1-4z) (1-2z)^{10} (1-z)^8 } \\
 &= 
\frac{(1-2z)^2}{(1-z)(1-4z)} \cdot {\left(\frac{1-z}{1-2z}\right)}^{10} 
\cdot \frac{1}{(1-z)^{17}},\label{zeta}
\end{align}
\smallskip
where in the second line we have written the zeta function as the
product of the zeta functions of one 2d solenoid $\mathbb{S}_2^2$,
ten 1d solenoids $\mathbb{S}_2$, and 17 extra fixed points. 

How can this be interpreted? A Robinson tiling generically consists
of a single, infinite order supertile. Such tilings project 1-to-1
to the solenoid $\mathbb{S}_2^2$. However, a Robinson tiling can 
consist also of 
two infinite order supertiles, which are separated by a horizontal 
or vertical row of tiles without any crosses. These are the tilings 
where the projection to $\mathbb{S}_2^2$ is not 1-to-1. A separating 
row of tiles can be decorated with a single blue line, or a double 
line with the second line (red) on either side of the middle blue line,
and all three 
cases can be combined with arrows in one or the other direction. 
All six possibilities, everything else being the same, project to 
the same point on $\mathbb{S}_2^2$. Moreover, if we take the 
translation orbit along the defect line, we obtain a whole 1d 
sub-solenoid $\mathbb{S}_2$ of such 6-tuples. So, in addition to 
the 2d solenoid $\mathbb{S}_2^2$, the hull contains 5 extra 1d 
solenoids $\mathbb{S}_2$in horizontal and 5 in vertical direction. 
Further, there are 28 fixed points of the substitution, 
consisting of 4 infinite order supertiles, which all 
project to the origin of $\mathbb{S}^2_2$. The 2d solenoid and the 
10 extra 1d solenoids contain one such fixed point each, so that in 
addition to those there must be 17 further ones, which all show up 
in the zeta function (\ref{zeta}). We finally note that the 
structure of the hull is in line with the interpretation of 
\cite{BS11}, were terms $\mathbb{Z}[{\scriptstyle\frac12}]$ in 
$H^2$ are associated with extra 1d sub-solenoids $\mathbb{S}_2$ 
in the hull.

\begin{acknowledgement}
The author would like to thank J.M.~Taylor for sharing her ideas
on the overlapping substitution rules for the Robinson tilings.
\end{acknowledgement}
%

%
%
%

\end{document}